\documentclass{amsart}

\usepackage{amsmath,amssymb,amscd,amsfonts}
\usepackage[all]{xy}

\newtheorem{theorem}{Theorem}
\newcommand{\bt}{\begin{theorem}}
\newcommand{\et}{\end{theorem}}
\newtheorem{lemma}{Lemma}
\newcommand{\bl}{\begin{lemma}}
\newcommand{\el}{\end{lemma}}
\newtheorem{corollary}{Corollary}
\newcommand{\bc}{\begin{corollary}}
\newcommand{\ec}{\end{corollary}}
\newcommand{\beq}{\begin{equation}}
\newcommand{\eeq}{\end{equation}}
\newcommand{\benum}{\begin{enumerate}}
\newcommand{\eenum}{\end{enumerate}}

\newcommand{\Z}{\ensuremath{\mathbf Z}}

\DeclareMathOperator{\qqand}{\qquad\text{and}\qquad}
\DeclareMathOperator{\kernel}{\text{kernel}}

\title[Not all groups are LEF groups]{Not all groups are LEF groups, or \\
can you know if a group is infinite?}

\author{Melvyn B. Nathanson}
\address{Lehman College (CUNY), Bronx, New York 10468}
\email{melvyn.nathanson@lehman.cuny.edu}

\date{\today}


\begin{document}

\maketitle

\begin{abstract}
This is an introduction to the class of groups 
that are locally embeddable into finite groups.
\end{abstract}

\section{Finite or infinite?}

A simple question:  Do the finite subsets of a group tell us if the group is infinite?  
Assume that we can only see the finite subsets of a group, 
and, also, that we can determine if a finite subset is a subset of some finite group.  
This means that we can answer the following question.    
Let $A$ be a finite subset of a group $G$.  
Does there exist a finite group $H$ and a \emph{partial homomorphism}\index{partial homomorphism} 
$f:A\rightarrow H$ that is one-to-one.  
A partial homomorphism from a subset $A$ of a group to a group $H$ is 
a function $f:A\rightarrow H$ such that, 
if $a,b \in A$ and $ab \in A$, then $f(ab) = f(a)f(b)$.
A one-to-one partial homomorphism is also called a \index{local embedding}\emph{local embedding}.
Of course, if the group $G$ is finite, then, for every subset $A$ of $G$, the restriction of the identity 
homomorphism on $G$ to the subset $A$ is a local embedding into a finite group. 

Does there exist an infinite group $G$ such that every finite subset of $G$ looks like 
(equivalently, can be partially embedded into) a subset of a finite group?  
Does there exist an infinite group $G$ in which some finite subset of $G$ is not also a subset of a finite group?  

Theorem~\ref{LEF:theorem:Thompson} answers the second question.
The following example  answers the first question.  
Let $A$ be a nonempty finite subset of the infinite abelian group \Z. 
Choose an integer 
\[
m > \max\{ |a-b|: a,b \in A\} = \max(A) - \min(A).
\] 
Consider the function $f:A \rightarrow \Z/m\Z$ defined by $f(a) =a + m\Z$ for all $a \in A$.  
This is a partial homomorphism because it is the restriction of the canonical homomorphism 
$a\mapsto a+m\Z$ from \Z\ to $\Z/m\Z$.  
For $a,b \in A$, we have $f(a) = f(b)$ if and only if $a\equiv b \pmod{m}$ if and only if 
$m$ divides $|a-b|$.   The inequality $|a-b| < m$ implies 
that $f(a) = f(b)$ if and only if $a = b$, and so $f$ is a local embedding.  
Thus, every finite subset of the infinite group $\Z$ can be embedded into a finite cyclic group.  
By looking only at finite subsets, we cannot decide if \Z\ is infinite.  

Let us call a group $G$ \emph{locally embeddable into finite groups}, or an \emph{LEF group}, 
if every finite subset of $G$ can be embedded into a finite group.  
Mal'cev~\cite{malc40} introduced this concept in general algebraic structures.   
Vershik and Gordon~\cite{vers-gord97} extended it to groups, and obtained many 
fundamental results.   

Here are two classes of LEF groups.  

\bt
Every locally finite group is an LEF group.  
Every abelian group is an LEF group.
\et

\begin{proof}
A group is \emph{locally finite} if every finite subset generates a finite group. 
For such groups, the proof is immediate from the definition.  

For abelian groups, the proof  follows easily from the structure theorem for finitely generated abelian groups, 
and an easy modification of the preceding argument that \Z\ is an LEF  group.
\end{proof}

It is natural to ask:  Is every infinite group an LEF group, or does there exist an infinite group 
that is not an LEF group?

\section{Finitely presented groups} 
Let $W$ be a group with identity $e$, and let $X$ be a subset of $W$ that generates $W$.  
We assume that $e \notin X$.  
The \emph{length} of an element  $w \in W$ with $w \neq e$ 
is the smallest  positive integer $k = \ell(w)$ such that there is a representation of $w$ in the form 
\beq           \label{LEF:w1}
w = x_1^{\varepsilon_1} x_2^{\varepsilon_2}  \cdots x_k^{\varepsilon_k} 
\eeq
where 
\beq                      \label{LEF:w2}
x_i \in X \text{ and } \varepsilon_i \in \{1,-1\} \text{ for } i=1,\ldots, k.
\eeq
We define $\ell(e) = 0$.   
Note that  $\ell(w) = 1$ if and only if $w = x$ or $w = x^{-1}$ for some $x \in  X$.  

For every nonnegative integer $L$, we define the ``closed ball'' 
\[
B_{L} = \{w \in W: \ell(w) \leq L \}.
\]
We have
\[
B_0 = \{e\} \qqand B_1 = \{e\} \cup \left\{x^{\varepsilon}: x \in X \text{ and } \varepsilon \in \{1,-1\} \right\}.
\]
If  the generating set $X$ is finite, then, for every $L$, the group $W$ contains 
only finitely many elements $w$ of length $\ell(w) \leq L$, and so 
the $B_L$ is a finite subset of $W$.  

If $w \in B_L$, then $w$ satisfies~\eqref{LEF:w1} and~\eqref{LEF:w2} for some $k \leq L$.   
For all $j = 1,\ldots, k$,  the partial product 
\[
w_j = x_1^{\varepsilon_1} x_2^{\varepsilon_2}  \cdots x_j^{\varepsilon_j}  
\]
has length $\ell(w_j) \leq j \leq L$, and so $w_j \in B_L$.
(We observe that if $\ell(w_j) < j$, then $\ell(w) < k$, which is absurd.  
Therefore, $\ell(w_j) = j$ for all $j \in \{1,\ldots, k\}$.) 
Let $w_0 = e$.  Note that $w = w_k$ and that 
\[
w_j = w_{j-1} x_j^{\varepsilon_j} 
\]
for all $j \in \{ 1,\ldots, k\}$.  
If $f:B_L \rightarrow H$ is a partial homomorphism, then
\begin{align*}
f(w) & = f(w_{k-1} x_k^{\varepsilon_k} ) = f(w_{k-1})f( x_k^{\varepsilon_k} ) \\
& = f(w_{k-2} x_{k-1}^{\varepsilon_{k-1}})  f( x_k^{\varepsilon_k} )
= f(w_{k-2} ) f\left( x_{k-1}^{\varepsilon_{k-1}} \right)  f( x_k^{\varepsilon_k} ) \\
& = \cdots \\
& =  f\left( x_1^{\varepsilon_1} \right)  f( x_2^{\varepsilon_2} ) \cdots  f\left( x_{k-1}^{\varepsilon_{k-1}} \right)  f( x_k^{\varepsilon_k} )
\end{align*}
For partial products in finite groups, see Nathanson~\cite{nath76b}.

Let $X$ be a nonempty  set, and  let $F(X)$ be the free group generated by $X$.
Let $R$ be a nonempty subset of $F(X)$. 
The \emph{normal closure}\index{normal closure} of $R$ in $F(X)$, denoted $N(R)$, 
is the smallest normal subgroup of $F(X)$ that contains $R$.  
The subgroup $N(R)$ is generated by the set 
\[
\left\{wr^{\varepsilon}w^{-1} : w \in F(X), \  r \in R,  \ \varepsilon \in \{1,-1\} \right\}.
\]
A group $G$ is \emph{finitely presented}\index{finitely presented group}\index{group!finitely presented} 
if 
\[
G =  \langle X;R \rangle = F(X)/N(R)
\]
where $F(X)$ is the free group generated by a finite set $X$ and the subgroup $N(R)$ 
is the normal closure of a finite subset $R$ of $F(X)$.  
If $\pi:F(X) \rightarrow G$ is the canonical homomorphism, then the set 
\[
X^{\ast} = \pi(X) = \{ xN(R):x\in X\}
\]
generates $G$.

The following result is Proposition 1.10 in Pestov and Kwiatkowska~\cite{pest-kwia13}.

\bt                     \label{LEF:theorem:PK}
Let $G$ be a finitely presented infinite group.  
If $G$ is an LEF group, then $G$ contains a nontrivial proper normal subgroup.  
Equivalently, a finitely presented  infinite simple group is not an LEF group.  
\et

\begin{proof}
Let $G = \langle X;R \rangle = F(X)/N$ be a finitely presented infinite group, 
where $F(X)$ is the free group generated by a finite set $X$, and $N = N(R)$ is the normal closure of 
a finite subset $R$ of $F(X)$.   
Let $e_F$ be the identity in $F(X)$.   The identity in $G$ is $e_G = e_FN = N$.  
The canonical homomorphism $\pi:F(X) \rightarrow G$ 
is defined by $\pi(w) = wN$ for all $w \in F(X)$.  

Choose an integer $L$ such that 
\[
L \geq \max\{ \ell(w): w \in X \cup R\}.
\]
The closed ball 
\[
B_{L} = \{w \in F(X): \ell(w) \leq L \}
\]
is a finite subset of $F(X)$.
We have   
\[
\{e_F\} \cup X \cup X^{-1} \cup R  \subseteq B_{L}.
\] 
The set  
\[
A= \pi\left(  B_{L} \right) \subseteq G 
\]
is a finite subset of $G$ that contains $X^* = \pi(X)$.   
Also,  $e_G = \pi(e_F) = N \in A$.

If $G$ is an LEF group, then there exist a finite group $H$ 
and a local embedding $f$ of  $A$ into $H$.  
Let $e_H$ be the identity in $H$.  
For all $x \in X$, we have $\pi(x) \in A$  and so 
\[
f \pi(x) \in H.
\]
By the universal property of a free group, there exists a unique homomorphism 
\[
f^*: F(X) \rightarrow H
\]
such that 
\[
f^*(x) = f\pi(x)
\]
for all $x \in X$.   
The subgroup 
\[
N^* = \kernel(f^*)  
\]
is a normal subgroup of $F(X)$.
We shall prove that 
\beq      \label{LEF:N-N^*}
N \subsetneq N^* \subsetneq F(X).
\eeq
The diagram is 
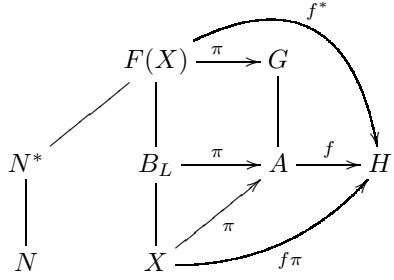
\[
\xymatrix{
& F(X) \ar[r]^{\pi}  \ar@{-}[d]    \ar@/^3pc/[rrd]^{f^*}   & G  \ar@{-}[d] & \\       
N^* \ar@{-}[ru] \ar@{-}[d]  & B_L \ar[r]^{\pi}   \ar@{-}[d]  & A \ar[r]^{f} & H \\
N & X \ar[ru]_{\pi}     \ar@/_1pc/[rru]_{f \pi}   & & 
}
\]

If $N^* = F(X)$, then $x \in  N^*$ for all $x \in X$.  
Because $X \subseteq  B_L$ and $\pi(x) = xN \in A$, we have 
\[
f(xN) = f\pi(x) = f^*(x) = e_H =  f(N).  
\]
Because $f$ is one-to-one and $f(xN) =  f(N)$,
it follows that $\pi(x) = xN = N$ for all $x \in X$.
The set $\pi(X)$ generates $G$, and so $G = \{ N \}$ is the trivial group, 
which is absurd.
Therefore, $N^*$ is a proper normal subgroup of $F(X)$.

Next we prove that $N^*$ contains $N$.  
Let $r \in R$.  There is a nonnegative integer $k = \ell(r)  \leq L$  
such that 
\[
r = \prod_{i=1}^k x_i^{\varepsilon_i}
\]
where $(x_i)_{i=1}^k$ is a sequence of elements 
of $X$ and $(\varepsilon_i)_{i=1}^k$ is a sequence of elements of $\{1,-1\}$.  

Because $r \in R \subseteq N$, we have  $rN = N$ and 
\begin{align*}
f^*(r) & = f^*\left( \prod_{i=1}^k x_i^{\varepsilon_i} \right) 
 = \prod_{i=1}^k f^*\left(  x_i\right)^{\varepsilon_i} 
  = \prod_{i=1}^k f \pi \left(  x_i \right)^{\varepsilon_i} \\
&    = \prod_{i=1}^k f \left(  x_i N\right)^{\varepsilon_i} 
 = f \left(  \prod_{i=1}^k  x_i^{\varepsilon_i} N\right) 
 = f \left( rN \right) \\
& = f \left( N \right) = e_H.
\end{align*}
Therefore, $r \in N^*$.   Because $R \subseteq N^*$   
and $N^*$ is a normal subgroup of $F(X)$, 
it follows that $N^*$ contains $N$, which is the 
normal closure of $R$, and so $N \subseteq N^*$.  

Finally, if $N = N^* = \kernel(f^*)$, then $G = F(X)/N = F(X)/N^*$ is isomorphic 
to a subgroup of the finite group $H$,
and so $G$ is finite, which is absurd. 
Therefore, $N$ is a proper subgroup of $N^*$.

This proves relation~\eqref{LEF:N-N^*}.  The correspondence theorem in group theory implies  
that $N^*/N$ is a nontrivial proper normal subgroup of $G$, 
and so $G$ is not a simple group.  
It follows that no  finitely presented infinite simple group is an LEF group.
This completes the proof.  
\end{proof}

\bt                  \label{LEF:theorem:Thompson}
There exist infinite groups that are not LEF groups.  
In particular, the Thompson groups $T$ and $V$ are not  LEF groups.
\et

\begin{proof}
The Thompson groups $T$ and $V$ are finitely presented  infinite simple  groups
(Cannon, Floyd, and Parry~\cite{cann-floy-parr96},Cannon and Floyd~\cite{cann-floy11}).  
\end{proof}

For recent work, including other examples of groups that are not LEF groups, 
see~\cite{corn-guyo-pits07} 
and~\cite{corn13}.

\def\cprime{$'$} \def\cprime{$'$} \def\cprime{$'$}
\providecommand{\bysame}{\leavevmode\hbox to3em{\hrulefill}\thinspace}
\providecommand{\MR}{\relax\ifhmode\unskip\space\fi MR }
\providecommand{\MRhref}[2]{%
  \href{http://www.ams.org/mathscinet-getitem?mr=#1}{#2}
}
\providecommand{\href}[2]{#2}

\end{document}